# Uniform saddlepoint approximations for ratios of quadratic forms

RONALD W. BUTLER[1] and MARC S. PAOLELLA[2]

[1]*Department of Statistical Science, Southern Methodist University, Dallas, TX 75275-0332, USA.* E-mail: *walrus@stat.colostate.edu*

[2]*Swiss Banking Institute, University of Zurich, CH-8032 Zurich, Switzerland.*
E-mail: *paolella@isb.unizh.ch*

Ratios of quadratic forms in correlated normal variables which introduce noncentrality into the quadratic forms are considered. The denominator is assumed to be positive (with probability 1). Various serial correlation estimates such as least-squares, Yule–Walker and Burg, as well as Durbin–Watson statistics, provide important examples of such ratios. The cumulative distribution function (c.d.f.) and density for such ratios admit saddlepoint approximations. These approximations are shown to preserve uniformity of relative error over the entire range of support. Furthermore, explicit values for the limiting relative errors at the extreme edges of support are derived.

## 1. Introduction

Consider the ratio of quadratic forms

$$R = \frac{\epsilon' \mathbf{A} \epsilon}{\epsilon' \mathbf{B} \epsilon}, \tag{1}$$

where, without loss in generality, $\mathbf{A}$ and $\mathbf{B}$ are assumed to be $n \times n$ symmetric. Let $\epsilon \sim \mathbf{N}(\mu, \mathbf{I}_n)$ and suppose that $\mathbf{B}$ is also positive semidefinite, thereby ensuring that the denominator is positive with probability one. There is no loss in generality in having the covariance of $\epsilon$ as the identity. This is because, if the distribution of $\epsilon$ were $\mathbf{N}(\mu, \Sigma)$, then (1) describes the model with $\Sigma^{1/2} \mathbf{A} \Sigma^{1/2}$ and $\Sigma^{1/2} \mathbf{B} \Sigma^{1/2}$ replacing $\mathbf{A}$ and $\mathbf{B}$, respectively, and $\Sigma^{-1/2} \mu$ replacing $\mu$ in the distributional assumption on $\epsilon$. Thus, model (1) incorporates all dependence among the components of $\epsilon$, as well as noncentrality that occurs when $\mu \neq \mathbf{0}$.

Various types of saddlepoint approximations for the distribution (c.d.f.) and density of $R$ have been proposed, beginning with the seminal work on serial correlations in Daniels (1956). Further marginal distributional approximations are given in McGregor (1960), Phillips (1978), Jensen (1988), Wang (1992), Lieberman (1994a, 1994b) and Marsh (1998). Joint distributional approximations for the set of serial correlations comprising the correlogram were initiated by Daniels (1956) and continued in Durbin (1980) and Butler and Paolella (1998).

The main contributions of the current paper are in establishing the uniformity of relative errors for the saddlepoint cdf and density approximations in the right tail when used with univariate ratios $R$ in a class $\mathcal{C}_\mathcal{R}$ that is defined below. This class encompasses all of the examples in the aforementioned papers. Expressions for the limiting relative error are given as the right edge of support for $R$ is approached and the sample size $n$ is held fixed. These expressions are explicit in the more elementary settings in which a certain defining eigenvalue is simple and mostly implicitly defined when it is multiple. The noncentral beta distribution provides an important example in which the limiting error is explicit, but the defining eigenvalue is multiple.

The left tail of $R$ is dealt with by changing $\mathbf{A}$ to $-\mathbf{A}$, thus switching the left tail of $R$ to the right tail of $-R$. The results for the right tail of $-R$ can now provide similar uniformity results for the left tail, when applicable. If $-R$ is a member of class $\mathcal{C}_\mathcal{R}$, then we say that $R$ is in $\mathcal{C}_\mathcal{L}$. However, for the most part, this paper concentrates on the class $\mathcal{C}_\mathcal{R}$.

The class of ratios $\mathcal{C}_\mathcal{R}$ is technically characterized in terms of a sequence of largest eigenvalues. Let $(\mathsf{l}, \mathsf{r})$ be the support of $R$ with $\mathsf{r}$ possibly infinite and define $\lambda_n(r)$ as the largest eigenvalue of $\mathbf{A} - r\mathbf{B}$ for $r \in (\mathsf{l}, \mathsf{r})$. The class $\mathcal{C}_\mathcal{R}$ is characterized as those ratios $R$ whose matrices $\mathbf{A}$ and $\mathbf{B}$ admit the limit

$$0 = \lim_{r \to \mathsf{r}} \lambda_n(r) \qquad (2)$$

with dimension $n$ fixed. The class $\mathcal{C}_\mathcal{R}$ contains the subset $\mathcal{B}$ that consists of all ratios with bounded support $(\mathsf{l}, \mathsf{r})$; this is a property of $R$ that is guaranteed when $\mathbf{B} > \mathbf{0}$, or positive definite. Such ratios include the Durbin–Watson statistics, as well as the Yule–Walker and Burg estimators of serial correlation with arbitrary lag computed from least-squares residuals. If $\mathbf{B} \geq \mathbf{0}$ has at least one zero eigenvalue, then $\mathsf{r}$ may be finite or infinite. The portion of $\mathcal{C}_\mathcal{R}$ with $\mathsf{r} = \infty$ includes least-squares estimators of serial correlation in various types of models with arbitrary lag and computed from residuals with trend or covariates removed. Such models include those with autoregressive lag in the dependent variable and those with lag in the additive noise.

Large *sample space* asymptotics of the type considered here have not been previously considered for the class $\mathcal{B}$. The only previous consideration for a member of the class $\mathcal{C}_\mathcal{R} - \mathcal{B}$ is in Jensen (1988), (1995), Chapter 9.4. Results obtained there are in agreement with those below for the least-squares estimator of lag-one serial correlation when the time series is a mean-zero $AR(1)$ model.

The class $\mathcal{C}_\mathcal{R}$ excludes $F$-statistic and Satterthwaite-type ratios which have been considered in Butler and Paolella (2002). In this work, $\lambda_n(r) > 0$ does not depend on $r$. For this setting, saddlepoint uniformity is also maintained; however, a different asymptotic saddlepoint behavior results from these different assumptions.

Some alternative large *sample size* asymptotics for the lag-one least-squares estimator, showing that the error is $O(n^{-1})$ and $O(n^{-3/2})$ on compact sets as $n \to \infty$, are given in Lieberman (1994b) and Jensen (1995), Chapter 9.4, respectively, when $\mu = \mathbf{0}$. Such asymptotics, in which $n \to \infty$, are not considered in this paper.

Since most of the proofs in this paper are long and technical, they have been relegated to the accompanying technical report, Butler and Paolella (2007).

## 2. Saddlepoint approximations

### 2.1. Distribution theory

The cdf for $R$ in the most general setting with noncentrality is

$$\Pr(R \leq r) = \Pr\left(\frac{\boldsymbol{\epsilon}'\mathbf{A}\boldsymbol{\epsilon}}{\boldsymbol{\epsilon}'\mathbf{B}\boldsymbol{\epsilon}} \leq r\right) = \Pr(\boldsymbol{\epsilon}'(\mathbf{A} - r\mathbf{B})\boldsymbol{\epsilon} \leq 0) \tag{3}$$
$$= \Pr(X_r \leq 0),$$

where $X_r = \boldsymbol{\epsilon}'(\mathbf{A} - r\mathbf{B})\boldsymbol{\epsilon}$

$$\mathbf{A} - r\mathbf{B} = \mathbf{P}_r' \boldsymbol{\Lambda}_r \mathbf{P}_r, \tag{4}$$

where $\mathbf{P}_r$ is orthogonal and $\boldsymbol{\Lambda}_r = \text{diag}(\lambda_1, \ldots, \lambda_n)$, with

$$\lambda_1 = \lambda_1(r) \leq \cdots \leq \lambda_n = \lambda_n(r)$$

consisting of the ordered eigenvalues of (4). Whenever convenient, we suppress the dependence of the various quantities on $r$. The distribution of $X_r$ is therefore

$$X_r = \sum_{i=1}^{n} \lambda_i \chi^2(1, \nu_i^2), \tag{5}$$

where $\{\nu_i^2\}$ are determined as $(\nu_1, \ldots, \nu_n)' = \boldsymbol{\nu}_r = \mathbf{P}_r \boldsymbol{\mu}$ and represent the noncentrality parameters of the independent noncentral $\chi_1^2$ variables specified in (5). The ordered values of $\{\lambda_i\}$ are in one-to-one correspondence with the components of $\boldsymbol{\nu}_r$ specified through the particular choice of $\mathbf{P}_r$. Notationally, we use $\chi_k^2$ for the central chi-square instead of $\chi^2(k, 0)$.

Before proceeding with the development of a saddlepoint approximation for the distribution of $R$, we must first characterize the support of $R$, its relationship to the eigenvalues $\lambda_1(r)$ and $\lambda_n(r)$ and the convergence strip for the moment generating function of $X_r$.

**Lemma 1.** *All of the eigenvalues of $\boldsymbol{\Lambda}_r$ are strictly decreasing in $r$ when $\mathbf{B} > \mathbf{0}$ and decreasing when $\mathbf{B} \geq \mathbf{0}$.*

**Lemma 2.** *The distribution of $R$ is degenerate at a single point if and only if $\mathbf{A} = c\mathbf{B}$ for some scalar constant $c$.*

A description of the support of $R$ requires the consideration of the various cases involved which depend on eigenvalue decompositions of $\mathbf{A}$ and $\mathbf{B}$. Suppose that $\mathbf{B}$ has $p \geq 0$ zero eigenvalues and let $\mathbf{O}'_\mathbf{B}$ be the orthogonal matrix of eigenvectors for $\mathbf{B}$ such that

$$\mathbf{O_B B O'_B} = \begin{pmatrix} \mathbf{\Lambda_B} & \mathbf{0} \\ \mathbf{0} & \mathbf{0}_{p \times p} \end{pmatrix}.$$

Denote

$$\mathbf{O_B A O'_B} = \begin{pmatrix} \mathbf{C}_{11} & \mathbf{C}_{12} \\ \mathbf{C}_{21} & \mathbf{C}_{22} \end{pmatrix},$$

where $\mathbf{C}_{11}$ is $(n-p) \times (n-p)$ and $\mathbf{C}_{22}$ is $p \times p$. Let $N(\mathbf{C}_{12})$ denote the null space in $\Re^p$ for matrix $\mathbf{C}_{12}$.

**Lemma 3.** *The support of $R$ is specified in the following set of cases.*

1. *Suppose that $\mathbf{B} > 0$, hence $p = 0$, and that $\mathbf{A}$ has rank of at least one. The support of $R$ is then the finite interval $(\mathsf{l}, \mathsf{r})$ with $\mathsf{l}$ and $\mathsf{r}$ being the smallest and largest eigenvalues of $\mathbf{B}^{-1}\mathbf{A}$.*
2. *If $p \geq 1$, so that $\mathbf{B}$ has at least one zero eigenvalue, then the right edge $\mathsf{r}$ is given as follows:*
   (a) *if $\mathbf{C}_{22}$ has a positive eigenvalue, then $\mathsf{r} = \infty$;*
   (b) *if $\mathbf{C}_{22} < \mathbf{0}$ then $\mathsf{r} < \infty$ and $\mathsf{r}$ is the largest eigenvalue of $\mathbf{\Lambda}_\mathbf{B}^{-1}(\mathbf{C}_{11} - \mathbf{C}_{12}\mathbf{C}_{22}^{-1}\mathbf{C}_{21})$;*
   (c) *if $\mathbf{C}_{22} \leq \mathbf{0}$ and $\mathbf{C}_{22}$ has at least one zero eigenvalue, then $\mathsf{r} = \infty$ if $N(\mathbf{C}_{22}) \not\subseteq N(\mathbf{C}_{12})$; otherwise $\mathsf{r} < \infty$ and is the largest eigenvalue of $\mathbf{\Lambda}_\mathbf{B}^{-1}(\mathbf{C}_{11} - \mathbf{C}_{12}\mathbf{O}_{\mathbf{C}1}\mathbf{\Lambda}_\mathbf{C}^{-1}\mathbf{O}'_{\mathbf{C}1}\mathbf{C}_{21})$. Here, $\mathbf{O}'_\mathbf{C} = (\mathbf{O}_{\mathbf{C}1}, \mathbf{O}_{\mathbf{C}2})$ consists of the eigenvectors of $\mathbf{C}_{22}$,*

$$\mathbf{O_C C}_{22}\mathbf{O'_C} = \begin{pmatrix} \mathbf{\Lambda_C} & \mathbf{0} \\ \mathbf{0} & \mathbf{0}_{m \times m} \end{pmatrix},$$

   $\mathbf{\Lambda_C} < \mathbf{0}$, $m$ *is the multiplicity of the zero eigenvalue and the columns of $\mathbf{O}_{\mathbf{C}1}$ and $\mathbf{O}_{\mathbf{C}2}$ consist of eigenvectors with nonzero and zero eigenvalues, respectively.*

Some of the settings described in Lemma 3 concern ratios that are not in the class $\mathcal{C}_\mathcal{R}$.

**Lemma 4.** *When considering the right tail, matrices $\mathbf{A}$ and $\mathbf{B}$ admit a ratio $R$ in the class $\mathcal{C}_\mathcal{R}$ only for cases 1, 2(b) or 2(c). When considering both the left and right tails, then the class $\mathcal{C}_\mathcal{R} \cap \mathcal{C}_\mathcal{L}$ encompasses case 1 and the special setting of case 2(c) in which $\mathbf{C}_{22} = \mathbf{0}$.*

Lemmas 3 and 4 are most easily understood by means of some simple examples. Consider an $F_{1,1}$ distribution for $R$. Then,

$$\mathbf{A} - r\mathbf{B} = \begin{pmatrix} 1 & 0 \\ 0 & -r \end{pmatrix}$$

and $\lambda_1(r) = -r$ with $\lambda_2(r) \equiv 1$. Clearly, this is not in the class $\mathcal{C}_\mathcal{R}$ nor in $\mathcal{C}_\mathcal{L}$. Since $\mathbf{C}_{22} = 1$, a scalar, this is case 2(a).

Next, consider $n = 2$ and the least-squares estimate of a lag-1 serial correlation in the simplest setting with $R = \epsilon_1 \epsilon_2 / \epsilon_1^2 = \epsilon_2/\epsilon_1$. Note that this has the Cauchy distribution when $\mu = 0$ and the support is $(l, r) = (-\infty, \infty)$. To see that this ratio is in the classes $\mathcal{C}_\mathcal{R}$ and $\mathcal{C}_\mathcal{L}$, note that

$$\mathbf{A} - r\mathbf{B} = \begin{pmatrix} -r & 1/2 \\ 1/2 & 0 \end{pmatrix} \tag{6}$$

and that the limiting eigenvalues are

$$\lim_{r \to -\infty} (-r - \sqrt{r^2+1})/2 = 0 = \lim_{r \to \infty} (-r + \sqrt{r^2+1})/2. \tag{7}$$

The example illustrates a case 2(c) ratio in which $\mathbf{C}_{12} = 1/2$ and $\mathbf{C}_{22} = 0$ are scalars and $N(\mathbf{C}_{22}) \nsubseteq N(\mathbf{C}_{12})$. The same results hold more generally with least-squares estimates of serial correlation from regression residuals.

**Lemma 5.** *Suppose that $R$ has a nondegenerate distribution in the class $\mathcal{C}_\mathcal{R}$, as described in Lemma 4, $\mathbf{B} \geq 0$ and $\mathbf{A}$ has rank of at least one. The upper range of support $\mathsf{r} \leq \infty$ for $R$, as given in cases 1, 2(b) and 2(c) of Lemma 3, solves $\lambda_n(\mathsf{r}) = 0$. If $\mathsf{r}$ is an interior point of the support of $R$, then the moment generating function of $X_r$ is*

$$M_{X_r}(s) = \left(\prod_{i=1}^n (1 - 2s\lambda_i)^{-1/2}\right) \exp\left\{s \sum_{i=1}^n \frac{\lambda_i \nu_i^2}{1 - 2s\lambda_i}\right\}, \tag{8}$$

*convergent in the neighborhood of zero given by*

$$\frac{1}{2\lambda_1(r)} < s < \frac{1}{2\lambda_n(r)}. \tag{9}$$

### 2.2. C.d.f. saddlepoint approximation

The saddlepoint approximation is based on the cumulant generating function (c.g.f.) for $X_r$, given by $K_{X_r}(s) = \ln M_{X_r}(s)$. The saddlepoint $\hat{s}$ is the unique root of

$$0 = K'_{X_r}(\hat{s}) = \sum_{i=1}^n \left(\frac{\lambda_i}{1 - 2\hat{s}\lambda_i} + \frac{\lambda_i \nu_i^2}{(1 - 2\hat{s}\lambda_i)^2}\right) \tag{10}$$

in the range (9). The approximation of Lugannani and Rice (1980) to first order is

$$\widehat{\Pr}_1(R \leq r) = \begin{cases} \Phi(\hat{w}) + \phi(\hat{w})\{\hat{w}^{-1} - \hat{u}^{-1}\}, & \text{if } 0 \neq \mathcal{E}[X_r], \\ \dfrac{1}{2} + \dfrac{K'''_{X_r}(0)}{6\sqrt{2\pi} K''_{X_r}(0)^{3/2}}, & \text{if } 0 = \mathcal{E}[X_r], \end{cases} \tag{11}$$

where $\Phi(\cdot)$ and $\phi(\cdot)$ denote the distribution and density function of a standard normal random variable, respectively, and

$$\hat{w} = \text{sgn}(\hat{s})\sqrt{-2K_{X_r}(\hat{s})}, \qquad \hat{u} = \hat{s}\sqrt{K''_{X_r}(\hat{s})}. \tag{12}$$

A second-order c.d.f. approximation is given in Butler and Paolella (2007).

### 2.3. Density saddlepoint approximation

The saddlepoint density approximation for $f_R(r)$, the density of $R$ at $r$, is derived in Butler (2007), Chapter 12.1, or Butler and Paolella (2007) as

$$\hat{f}_R(r) = \frac{J_r(\hat{s})}{\sqrt{2\pi K''_{X_r}(\hat{s})}} M_{X_r}(\hat{s}), \tag{13}$$

where $\hat{s}$ is the same saddlepoint used in the c.d.f. approximation and which solves (10). The factor $J_r(\hat{s})$ is computed from

$$J_r(s) = \text{tr}(\mathbf{I} - 2s\mathbf{\Lambda}_r)^{-1}\mathbf{H}_r + \nu'_r(\mathbf{I} - 2s\mathbf{\Lambda}_r)^{-1}\mathbf{H}_r(\mathbf{I} - 2s\mathbf{\Lambda}_r)^{-1}\nu_r \tag{14}$$

with $\mathbf{H}_r = \mathbf{P}_r\mathbf{B}\mathbf{P}'_r$. A second-order saddlepoint density is given in Butler (2007), page 383, or Butler and Paolella (2007).

*Example 6.* For matrices **A** and **B** in which $R \sim \text{Beta}(m/2, (n-m)/2)$, the saddlepoint density in (13) is

$$\hat{f}_R(r) = \frac{B(m/2, (n-m)/2)}{\hat{B}(m/2, (n-m)/2)} f_R(r),$$

where $\hat{B}$ is Stirling's approximation for the Beta function $B$.

## 3. Uniformity of the approximations in $r$

The relative errors of Lugannani and Rice's approximation in (11) and the density approximation in (13) can be shown to be uniform over $[0, \mathsf{r})$ when in the class $\mathcal{C}_{\mathcal{R}}$. These results follow as consequences of deriving their finite limiting ratios as $r \to \mathsf{r}$. The limiting ratios are derived in Theorems 9, 14 and 15 below. Our approach to computing these limiting ratios follows that used in Jensen (1988), (1995), Chapter 9.4, and generalizes these results to accommodate both noncentrality and the special concerns involving multiple eigenvalues.

The nature of these asymptotics is dependent on the multiplicity of the eigenvalue $\lambda_n(\mathsf{r}) = 0$, denoted as $m$. As a simple eigenvalue with $m = 1$, the limiting ratios are derived in Theorem 9. This is a common setting encountered when dealing with serial correlations. With $m \geq 2$, however, the asymptotics are more difficult and such results are deferred to Theorem 14. Examples of the multiple eigenvalue setting are also common and include least-squares estimates and Yule–Walker estimates for lag-$l$ serial correlation with $l \geq 2$. One important multiple eigenvalue example is the noncentral beta distribution discussed in Section 4.

The case 2(a) setting is not in $\mathcal{C}_\mathcal{R}$; however, the relative error can still be shown to be uniform over $[0, \infty)$; see Butler and Paolella (2007).

### 3.1. Simple eigenvalue $\lambda_n(r) = 0$

Suppose that $r < \infty$, under the circumstances of cases 1, 2(b) or 2(c). We assume here that $\mathbf{A} - r\mathbf{B}$ has a simple zero eigenvalue with multiplicity $m = 1$. For general $\mathbf{A}$ and $\mathbf{B}$, this multiplicity is often difficult to anticipate. Define

$$\nu_0 = \nu_n(r) := \mathbf{p}_n(r)'\boldsymbol{\mu}, \tag{15}$$

where $\mathbf{p}_n(r)$ is the eigenvector associated with the zero eigenvalue of $\mathbf{A} - r\mathbf{B}$.

The situation with $r = \infty$ is more complicated.

**Lemma 7.** *Suppose case* 2(c) *with* $r = \infty$. *Then,* $m$, *the multiplicity of zero eigenvalues in* $\{\lambda_i(\infty)\}$, *is the number of zero eigenvalues for* $\mathbf{C}_{22}$. *If* $m = 1$, *then*

$$\nu_0 = \nu_n(\infty) := \mathbf{o}_n' \mathbf{O}_{\mathbf{B}2}' \boldsymbol{\mu}, \tag{16}$$

*where* $\mathbf{o}_n$ *is the* $p \times 1$ *eigenvector associated with the zero eigenvalue of* $\mathbf{C}_{22}$, $\mathbf{O}_\mathbf{B}' = (\mathbf{O}_{\mathbf{B}1}, \mathbf{O}_{\mathbf{B}2})$ *and* $\mathbf{O}_{\mathbf{B}2}$ *is* $n \times p$ *and the orthonormal basis for the null space of* $\mathbf{B}$ *used to determine* $\mathbf{C}_{22}$.

The $AR(1)$ example in (6) with $n = 2$ provides a simple example. Here, $\mathbf{C}_{22}$ is the scalar 0 so that $\mathbf{o}_n' = 1$, and $\mathbf{O}_{\mathbf{B}2}' = (0, 1)$. Hence, $\nu_0 = \mu_2$.

**Lemma 8.** *Suppose that the conditions of Lemma 5 hold and let* $m = 1$. *Then, as* $r \to \mathsf{r} \leq \infty$, $\epsilon = \lambda_n(r) \to \lambda_n(\mathsf{r}) = 0$ *and* $\hat{s} = t_0/\epsilon + O(1) \to \infty$, *where*

$$t_0 = \frac{1}{4n}\{2n - 1 + \nu_0^2 - \sqrt{(\nu_0^2 + 2n - 1)^2 - (2n-1)^2 + 1}\} \tag{17}$$

*and* $\nu_0$ *is defined in* (15) *or* (16). *In addition,*

$$\hat{u} \to u_0 = \sqrt{\frac{n-1}{2} + \frac{2t_0^2}{(1-2t_0)^2} + \frac{4\nu_0^2 t_0^2}{(1-2t_0)^3}}. \tag{18}$$

**Theorem 9.** *Suppose that* $n \geq 2$, $R$ *has a nondegenerate distribution in* $\mathcal{C}_\mathcal{R}$, $\mathbf{B} \geq 0$ *and* $\mathbf{A}$ *has rank of at least one. If* $m = 1$, *then the limiting ratio of the true tail probability to* *its first order Lugannani–Rice approximation in* (11) *is*

$$\lim_{r \to \mathsf{r}} \frac{\Pr(R > r)}{\widehat{\Pr}_1(R > r)} = \frac{\sqrt{2\pi(1-2t_0)}(2t_0)^{(n-1)/2}u_0 e^{-\eta_2}}{B(1/2,(n+1)/2)(n/2)} {}_1F_1\left(\frac{n}{2};\frac{1}{2};\frac{\nu_0^2}{2}\right), \quad (19)$$

*where*

$$\eta_2 = \frac{\nu_0^2}{2(1-2t_0)} \quad (20)$$

*and parameters $t_0, u_0$ and $\nu_0$ are specified in* (17), (18) *and* (15), (16). *The first-order saddlepoint density has the same relative limit. All of these parameters are determined by $\nu_0$, so the right-hand side of* (19) *is a function of $\nu_0$ alone.*

In the central case with $\nu = 0$, the limiting ratio of tail probabilities in Theorem 9 is $\hat{B}(\frac{1}{2}, \frac{n-1}{2})/B(\frac{1}{2}, \frac{n-1}{2})$, where $\hat{B}$ is Stirling's approximation. This same limiting error was derived in Jensen (1995), Chapter 9.4, which considered the tail ratio for the distribution of the least-squares estimate in a mean zero $AR(1)$ model. Jensen's (9.4.7) is this value when the difference in notation is accounted for (our $n$ being $n+1$ in that paper).

As $\nu_0^2 \to \infty$, the limiting ratio in Theorem 9 is $\hat{\Gamma}(\frac{1}{2}, \frac{n-1}{2})/\Gamma(\frac{1}{2}, \frac{n-1}{2})\{1 + O(\nu_0^{-2})\}$, where $\hat{\Gamma}$ is Stirling's approximation. This follows from the large argument asymptotics for ${}_1F_1$ as given in 13.1.4 of Abramowitz and Stegun (1972).

Relative errors for the second order cdf and density approximations are also uniform in the right tail; see Butler and Paolella (2007).

## 3.2. Multiple eigenvalue $\lambda_n(\mathsf{r}) = 0$

In this setting, the asymptotics depend on the relative rates of convergence to zero for the multiple eigenvalues of $\mathbf{A} - r\mathbf{B}$ that approach 0 as $r \to \mathsf{r}$. If $m$ denotes its multiplicity, then $m \geq 1$, by the definition of $\mathsf{r}$. The allowable values of $m$ are $1 \leq m \leq n-1$, but not $m = n$. This latter value would make $\lim_{r \to \mathsf{r}}(\mathbf{A} - r\mathbf{B}) \equiv 0$, in which case the distribution of $R$ approaches a degenerate distribution, by Lemma 2. For unbounded ratios in case 2(c), the value of $m$ is the dimension of the null space for $\mathbf{C}_{22}$, whereas, for ratios in $\mathcal{B}$, the value of $m$ is less transparent.

We must first determine the relative rates at which the $m$ largest eigenvalues of $\mathbf{A} - r\mathbf{B}$ vanish as $r \to \mathsf{r}$ in the two separate settings, $\mathsf{r} < \infty$ and $\mathsf{r} = \infty$. For the former setting, general formulae for these relative rates are given in the next lemma. When $\mathsf{r} = \infty$, the relative rates must be determined on a case-by-case basis.

**Lemma 10.** *Suppose that $\mathsf{r} < \infty$ and 0 is an eigenvalue of multiplicity $m$ for $\mathbf{A} - \mathsf{r}\mathbf{B}$. Let the columns of the $n \times m$ matrix $\mathbf{U}_0$ be an orthonormal basis for the null space of $\mathbf{A} - \mathsf{r}\mathbf{B}$. Furthermore, denote the ordered eigenvalues of $\mathbf{U}_0'\mathbf{B}\mathbf{U}_0$ as $0 \leq \tau_{n-m+1} \leq \cdots \leq \tau_n$. If $\tau_n > 0$, then the limiting relative rates of convergence to zero for the $m$ largest eigenvalues* *of* $\mathbf{A} - r\mathbf{B}$ *are*

$$\lim_{r \to \mathsf{r}} \frac{\lambda_i(r)}{\lambda_n(r)} = \frac{\tau_i}{\tau_n} = \omega_i \qquad (21)$$

*for* $i = n - m + 1, \ldots, n$, *where* $0 \leq \omega_{n-m+1} \leq \cdots \leq \omega_n = 1$.

For the most common case, in which $\mathbf{B} > \mathbf{0}$, we have $\tau_{n-m+1} > 0$ so that $\omega_{n-m+1} > 0$. To deal with the $\mathsf{r} = \infty$ setting of case 2(c), we reparametrize

$$\mathbf{D}(\varepsilon) = (\mathbf{A} - r\mathbf{B})/r = \varepsilon \mathbf{A} - \mathbf{B} \qquad (22)$$

and let $\varepsilon = 1/r \to 0$. If $\lambda_i(\varepsilon)$ are the ordered eigenvalues of $\mathbf{A} - \varepsilon^{-1}\mathbf{B}$, then

$$\psi_i(\varepsilon) = \varepsilon \lambda_i(\varepsilon) \qquad (23)$$

are the ordered eigenvalues of (22).

**Lemma 11.** *Consider case* 2(c), *in which* $\mathsf{r} = \infty$, *and assume that the zero eigenvalue has multiplicity* $m$. $\lambda_{n-m+1}(r), \ldots, \lambda_n(r)$ *are then analytic at* $r = \infty$. *If* $\lambda'_n(\infty) > 0$, *then the relative rates of convergence are*

$$\omega_i = \frac{\partial \lambda_i(\varepsilon)/\partial \varepsilon|_{\varepsilon=0}}{\partial \lambda_n(\varepsilon)/\partial \varepsilon|_{\varepsilon=0}} = \frac{\partial^2 \psi_i(\varepsilon)/\partial \varepsilon^2|_{\varepsilon=0}}{\partial^2 \psi_n(\varepsilon)/\partial \varepsilon^2|_{\varepsilon=0}} \qquad (24)$$

*for* $i = n - m + 1, \ldots, n$.

The limiting noncentrality parameters $\{\nu_{0i} : i = n - m + 1, \ldots, n\}$ are more difficult to determine for $m \geq 2$ because they are expressed in terms of the limiting eigenvectors associated with the eigenvalues that vanish. In the case $\mathsf{r} < \infty$, it is intuitively clear and Lancaster (1964) has shown formally that these eigenvectors are smoothly defined as $r \to \mathsf{r}$. Let $\mathbf{P}_{2r}$ be $n \times m$ and consist of the last $m$ columns of $\mathbf{P}'_r$ which are the eigenvectors for the $m$ largest eigenvalues of $(\mathbf{A} - r\mathbf{B})$ (which increase in size with column number). $\mathbf{P}_{2r}$ is then continuous at $r = \mathsf{r}$ and the limiting noncentrality parameters are

$$(\nu_{0,n-m+1}, \ldots, \nu_{0,n})' = \mathbf{P}'_{2\mathsf{r}} \boldsymbol{\mu}. \qquad (25)$$

In the unbounded setting with $\mathsf{r} = \infty$, let the $n \times m$ matrix $\mathbf{P}_{2\varepsilon}$ consist of the eigenvectors corresponding to the largest $m$ eigenvalues of $\mathbf{D}(\varepsilon)$ in (22). The limiting noncentrality parameters are then given in (25) with $\mathbf{P}_{20} = \lim_{\varepsilon \to 0} \mathbf{P}_{2\varepsilon}$ replacing $\mathbf{P}_{2\mathsf{r}}$. In complicated practical examples where these computations are not explicit, these limiting eigenvectors are best computed numerically by using a small $\varepsilon > 0$.

*Example 12.* The least-squares estimate of a lag-2 serial correlation with $n = 3$ and zero mean has the form $R = \epsilon_1 \epsilon_3 / \epsilon_1^2$ and leads to the matrix

$$\mathbf{D}(\varepsilon) = \varepsilon \mathbf{A} - \mathbf{B} = \begin{pmatrix} -1 & 0 & \frac{1}{2}\varepsilon \\ 0 & 0 & 0 \\ \frac{1}{2}\varepsilon & 0 & 0 \end{pmatrix} = \mathbf{Q}'_\varepsilon \begin{pmatrix} \psi_-(\varepsilon) & 0 & 0 \\ 0 & 0 & 0 \\ 0 & 0 & \psi_+(\varepsilon) \end{pmatrix} \mathbf{Q}_\varepsilon,$$

where

$$\mathbf{Q}'_\varepsilon = \begin{pmatrix} 2\psi_-(\varepsilon)/\varepsilon & 0 & 2\psi_+(\varepsilon)/\varepsilon \\ 0 & 1 & 0 \\ 1 & 0 & 1 \end{pmatrix}, \qquad \psi_\pm(\varepsilon) = -\tfrac{1}{2} \pm \tfrac{1}{2}\sqrt{1+\epsilon^2}.$$

The eigenvectors in matrix $\mathbf{Q}'_\varepsilon$ have not been normalized as would be needed to use the notation $\mathbf{P}'_\varepsilon$. The limits of the eigenvalues are $\lim_{\varepsilon \to 0}\{\psi_-(\varepsilon), 0, \psi_+(\varepsilon)\} = (-1, 0, 0)$ and the limiting normed eigenvectors have $\mathbf{P}'_\varepsilon \to \mathbf{I}_3$ as $\varepsilon \to 0$. Note that $\partial \psi_+(\varepsilon)/\partial \varepsilon|_{\varepsilon=0} = 0$. Also, the eigenvalues of

$$\mathbf{C}_{22} = \mathbf{O}'_{\mathbf{B}2}\mathbf{A}\mathbf{O}_{\mathbf{B}2} = \begin{pmatrix} 0 & 1 & 0 \\ 0 & 0 & 1 \end{pmatrix} \begin{pmatrix} 0 & 0 & \tfrac{1}{2} \\ 0 & 0 & 0 \\ \tfrac{1}{2} & 0 & 0 \end{pmatrix} \begin{pmatrix} 0 & 0 \\ 1 & 0 \\ 0 & 1 \end{pmatrix} = \begin{pmatrix} 0 & 0 \\ 0 & 0 \end{pmatrix}$$

are both zero, as discussed in Lemma 11. The limiting rate

$$\omega_2 = \lim_{r \to \mathsf{r}} \frac{\lambda_2(r)}{\lambda_3(r)} = \lim_{\varepsilon \to 0} \frac{0}{\partial^2 \psi_+(\varepsilon)/\partial \varepsilon^2} = \frac{0}{1/2} = 0.$$

The limiting noncentrality parameters are

$$\begin{pmatrix} \nu_{02} \\ \nu_{03} \end{pmatrix} = \mathbf{P}'_{20}\mu = \begin{pmatrix} 0 & 1 & 0 \\ 0 & 0 & 1 \end{pmatrix} \mu = \begin{pmatrix} \mu_2 \\ \mu_3 \end{pmatrix}.$$

Expressions for the limiting relative errors depend on $\{\omega_i, \nu_{0i}\}$. All summations in the remainder of this subsection are over $\mathcal{S} = \{n-m+1, \ldots, n\}$.

**Lemma 13.** *Suppose that $R$ is in the class $\mathcal{C}_\mathcal{R}$ and let $m$ be the multiplicity of the zero eigenvalue of $\mathbf{A} - \mathsf{r}\mathbf{B}$. Then, as $r \to \mathsf{r}$, $\epsilon = \lambda_n(r) \to \lambda_n(\mathsf{r}) = 0$ and $\hat{s} = t_0/\epsilon + O(1) \to \infty$, where $t_0$ is the unique solution to*

$$0 = -\frac{n-m}{2t_0} + \sum_{i \in \mathcal{S}} \omega_i \left\{ \frac{1}{1-2t_0\omega_i} + \frac{\nu_{0i}^2}{(1-2t_0\omega_i)^2} \right\} \tag{26}$$

*in $(0, 1/2)$ with $\mathcal{S} = \{n-m+1, \ldots, n\}$. In addition,*

$$\hat{u} \to u_0 = \sqrt{\frac{n-m}{2} + 2t_0^2 \sum_{i \in \mathcal{S}} \omega_i^2 \left\{ \frac{1}{(1-2t_0\omega_i)^2} + \frac{2\nu_{0i}^2}{(1-2t_0\omega_i)^3} \right\}}. \tag{27}$$

**Theorem 14.** *Suppose that $n \geq 2$ and the conditions of Lemma 13 hold. Define the operator $D_0(X)$ to be the density of the random variable $X$ evaluated at zero. The limiting ratio of the true tail probability for $R$ to its first order Lugannani–Rice approximation in* (11) *is then*

$$\lim_{r \to \mathsf{r}} \frac{\Pr(R > r)}{\widehat{\Pr}_1(R > r)} = \sqrt{2\pi} D_0 \left\{ \sum_{i \in \mathcal{S}} \eta_{1i}\chi^2(1, 2\eta_{2i}) - \frac{1}{2u_0}\chi^2_{n-m+2} \right\}, \tag{28}$$

where the $\chi^2$ terms are independent random variables. Parameters $\eta_{1i}$ and $\eta_{2i}$, for $i \in \mathcal{S}$, are

$$\eta_{1i} = \frac{t_0 \omega_i}{u_0(1 - 2t_0\omega_i)}, \qquad \eta_{2i} = \frac{\nu_{0i}^2}{2(1 - 2t_0\omega_i)}.$$

The limiting ratio for the density approximation is

$$\lim_{r \to \mathsf{r}} \frac{f_R(r)}{\hat{f}_R(r)} = \frac{\sqrt{2\pi}}{W_J} \left[ \sum_{i \in \mathcal{S}} h_{ii} \eta_{3i} D_0 \left\{ \sum_{j \in \mathcal{S}} \eta_{1j} \chi^2(1, 2\eta_{2j}) + \eta_{1i}\chi_2^2 - \frac{1}{2u_0}\chi_{n-m}^2 \right\} \right. \tag{29}$$

$$\left. + \sum_{i \in \mathcal{S}} \sum_{j \in \mathcal{S}} \nu_{0i}\nu_{0j}\eta_{3i}\eta_{3j}h_{ij} D_0 \left\{ \sum_{k \in \mathcal{S}} \eta_{1k} \chi^2(1, 2\eta_{2k}) + \eta_{1i}\chi_2^2 + \eta_{1j}\chi_2^2 - \frac{1}{2u_0}\chi_{n-m}^2 \right\} \right],$$

where $\eta_{3i} = (1 - 2t_0\omega_i)^{-1}$,

$$W_J = \sum_{i \in \mathcal{S}} h_{ii} \eta_{3i} + \sum_{i \in \mathcal{S}} \sum_{j \in \mathcal{S}} \nu_{0i}\nu_{0j}\eta_{3i}\eta_{3j}h_{ij} > 0,$$

$(h_{ij}) = \mathbf{H}_\mathsf{r} = \lim_{r \to \mathsf{r}} \mathbf{P}_r \mathbf{B} \mathbf{P}_r$ and all $\chi^2$ variates in (29) are assumed to be independent.

In the case $m = 1$, the results of Theorem 14 reduce to those in Theorem 9.

## 4. Noncentral Beta $(\frac{m}{2}, \frac{n-m}{2})$ distribution

This distribution has

$$\mathbf{A} = \begin{pmatrix} \mathbf{I}_m & \mathbf{0} \\ \mathbf{0} & \mathbf{0} \end{pmatrix}$$

and $\mathbf{B} = \mathbf{I}_n$ so that $\mathsf{r} = 1$ and $\omega_i \equiv 1$. This leads to the explicit expression

$$t_0 = \frac{1}{2} + \frac{1}{4n}\{(\theta - m) - \sqrt{(\theta - m)^2 + 4\theta n}\},$$

where

$$\theta = \sum_{i \in \mathcal{S}} \nu_{0i}^2 = \sum_{i=1}^m \mu_i^2.$$

Furthermore,

$$u_0 = \sqrt{\frac{n-m}{2} + 2t_0^2 \left\{ \frac{m}{(1-2t_0)^2} + \frac{2\theta}{(1-2t_0)^3} \right\}}.$$

**Theorem 15.** *For a noncentral* $\text{Beta}(\frac{m}{2}, \frac{n-m}{2})$ *distribution with* $\min(m, n-m) \geq 1$, *the limiting ratio of the true tail probability to its first-order Lugannani–Rice approximation is*

$$RE = \frac{\sqrt{2\pi}(1-2t_0)^{m/2}(2t_0)^{(n-m)/2}u_0 e^{-\eta_2}}{B(m/2, (n-m)/2)(n-m)/2} {}_1F_1\left(\frac{n}{2}; \frac{m}{2}; \frac{\theta}{2}\right), \quad (30)$$

*where*

$$\eta_2 = \frac{\theta}{2(1-2t_0)}.$$

*The first-order saddlepoint density has the same relative error limit.*

In the central setting with $\theta = 0$, the value in (30) reduces to $\hat{B}(\frac{m}{2}, \frac{n-m}{2})/B(\frac{m}{2}, \frac{n-m}{2})$. This is consistent with the computation of the central $\text{Beta}(\frac{m}{2}, \frac{n-m}{2})$ density in Example 6.

As $\theta \to \infty$, the limiting ratio for (30) is $\hat{\Gamma}(\frac{1}{2}, \frac{n-m}{2})/\Gamma(\frac{1}{2}, \frac{n-m}{2})\{1 + O(\theta^{-1})\}$, which follows from the asymptotics for ${}_1F_1$ given in 13.1.4 of Abramowitz and Stegun (1972).

## 5. Examples

### 5.1. Serial correlations

Least-squares, Yule–Walker and Burg estimates for lag-$l$ correlations are considered in further detail in Butler and Paolella (2007) as members of the respective classes $\mathcal{C}_\mathcal{R}-\mathcal{B}$, $\mathcal{B}$ and $\mathcal{B}$. Both simple and multiple eigenvalue settings occur in such examples when there is no correction for mean.

In practice, serial correlations with arbitrary lag $l$ are generally computed from least-squares residuals and this often ensures that the largest eigenvalue of $\mathbf{A} - r\mathbf{B}$ has algebraic multiplicity one. Thus, the simpler situation for the large deviation errors occurs most often in practical data analysis.

### 5.2. Numerical example

Numerical confirmation of the large deviation errors in Theorem 9 is possible by considering the simplest model of Section 2.1. This is the least-squares estimate of lag one with $n = 2$ in a model without a location effect. Then, $R = \epsilon_2/\epsilon_1$ with $\epsilon_i \stackrel{\text{indep}}{\sim} N(\mu_i, 1)$. The exact density can be expressed as

$$\begin{aligned}
f_R(r) &= \frac{1}{2\pi} \int_{-\infty}^{\infty} |x| \exp\left\{-\frac{1}{2}(x-\mu_1)^2\right\} \exp\left\{-\frac{1}{2}(rx-\mu_2)^2\right\} dx \\
&= (\pi\delta)^{-1} \exp\left\{-\frac{1}{2}(\mu_1^2 + \mu_2^2)\right\} + \frac{\lambda\theta(\mu_1 + r\mu_2)}{\delta\sqrt{2\pi\delta}},
\end{aligned} \quad (31)$$

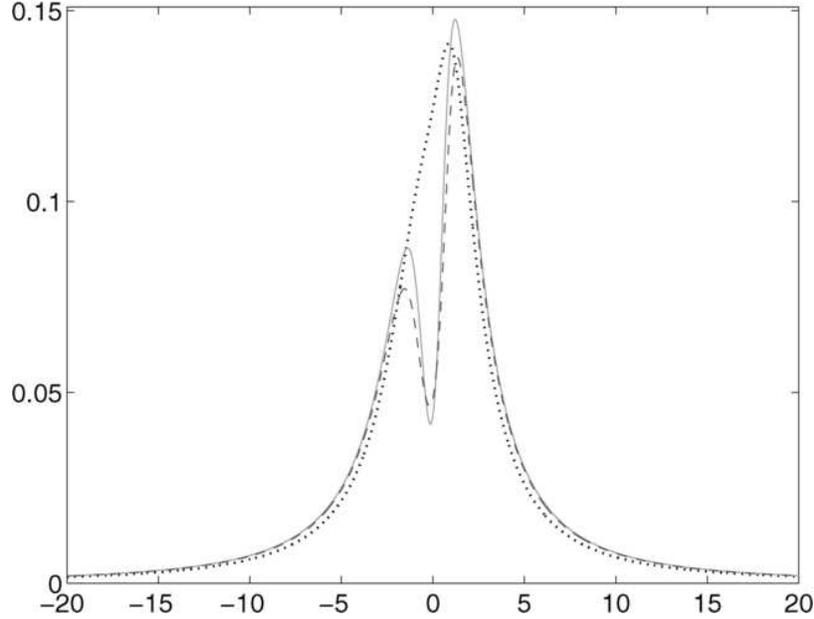

**Figure 1.** Exact density $f_R$ (solid), second order $\hat{f}_{R2}$ (dashed) and normalized $\bar{f}_R$ (dotted) approximations.

where

$$\delta = 1 + r^2, \qquad \theta = \operatorname{erf}\left(\frac{\mu_1 + r\mu_2}{\sqrt{2\delta}}\right), \qquad \lambda = \exp\left(-\frac{1}{2}\frac{(\mu_1 r - \mu_2)^2}{\delta}\right).$$

From Theorem 9, the limiting relative errors in the left and right tails are dependent on $\nu_0$ alone; in both tails, this value is $\nu_0 = \mu_2$ so that the limiting relative errors are the same in both tails, regardless of the values of $\mu_1$ and $\mu_2$.

The density (31) is both heavy-tailed and bimodal for $\mu_1 = 0.2$ and $\mu_2 = 2$. Figure 1 plots the exact density, the normalized version of $\hat{f}_R$ in (13) denoted by $\bar{f}_R$ and the second-order saddlepoint $\hat{f}_{R2}$ given in (20) of Butler and Paolella (2007). While both appear highly accurate in the tails, only the latter captures the bimodality. Figure 2 plots the ratio of the exact to the three approximate densities including $\hat{f}_R$, $\bar{f}_R$ and $\hat{f}_{R2}$. As $|r|$ increases, we have numerically confirmed that $f_R(r)/\hat{f}_R(r) \to 0.8222$, in agreement with the value computed via Theorem 9. This value is virtually achieved at $|r| = 10$. Both $\bar{f}_R$ and $\hat{f}_{R2}$ perform better than $\hat{f}_R$ in the tails, the latter most notably so.

The true cdf of $R$, or $F_R(r)$, must be computed from (31) using numerical integration. In this case, $|r|$ must be substantially larger before the same limiting ratio, as specified in Theorem 9, is reached. Figure 3 plots

$$\frac{F_R(r)}{\hat{F}_R(r)}1_{\{\hat{s}<0\}} + \frac{1-F_R(r)}{1-\hat{F}_R(r)}1_{\{\hat{s}>0\}} \quad \text{vs.} \quad r \tag{32}$$

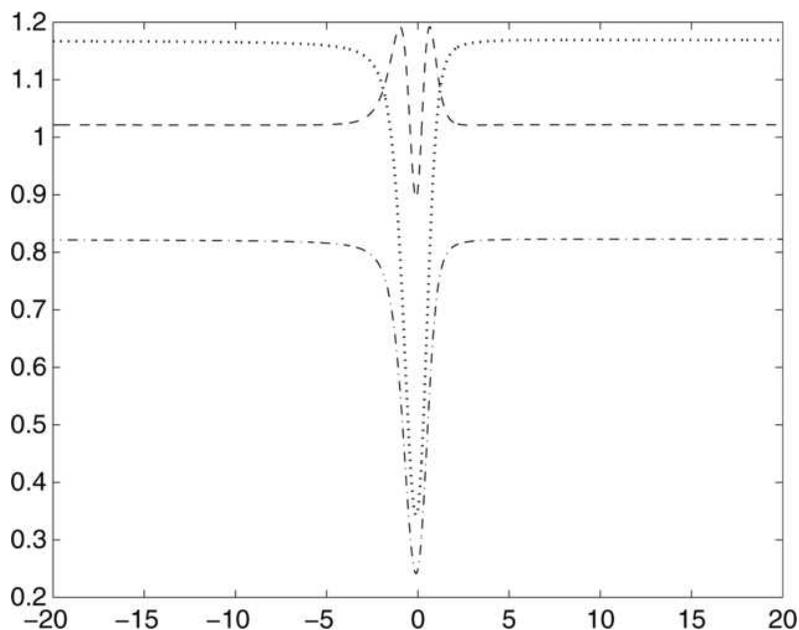

**Figure 2.** Error ratios $f_R/\hat{f}_{R2}$ (dashed), $f_R/\hat{f}_R$ (dashed-dot) and $f_R/\bar{f}_R$ (dotted).

with $\hat{F}_R(r)$ as the second-order approximation $\widehat{\Pr}_2$ in (17) of Butler and Paolella (2007) and also as $\widehat{\Pr}_1$ in (11). For these values of $\mu_i$, $\widehat{\Pr}_1$ is more accurate than $\widehat{\Pr}_2$ only in the range $-1.8 < r < 1.2$. At $r = -25,000$, $F_R(r)/\hat{F}_R(r) = 0.8226$ for $\widehat{\Pr}_1$, as given by Theorem 9, while, for $\widehat{\Pr}_2$, the ratio is 1.015. This latter ratio necessarily includes the factor $(1 + O_F)$, where $O_F$ approximates the limit of the second-order correction term.

If $\mu_1 = \mu_2 = 0$, then $R$ is Cauchy and the saddlepoint density reduces to $\hat{f}_R(r) = \sqrt{\pi/2} f_R(r)$. Thus, $\bar{f}_R$ is exact and the saddlepoint solution to $0 = K'_X(\hat{s})$ is given by $\hat{s} = r$. The relative error is, $\hat{B}(1/2, 1/2)/B(1/2, 1/2)$, in agreement with the large sample space theory.

## Acknowledgements

The research of R. Butler was supported by NSF Grants DMS-06-04318 and DMS-02-02284. The research of M. Paolella was partially carried out within the National Centre of Competence in Research "Financial Valuation and Risk Management" (NCCR FIN-RISK), a research program supported by the Swiss National Science Foundation.

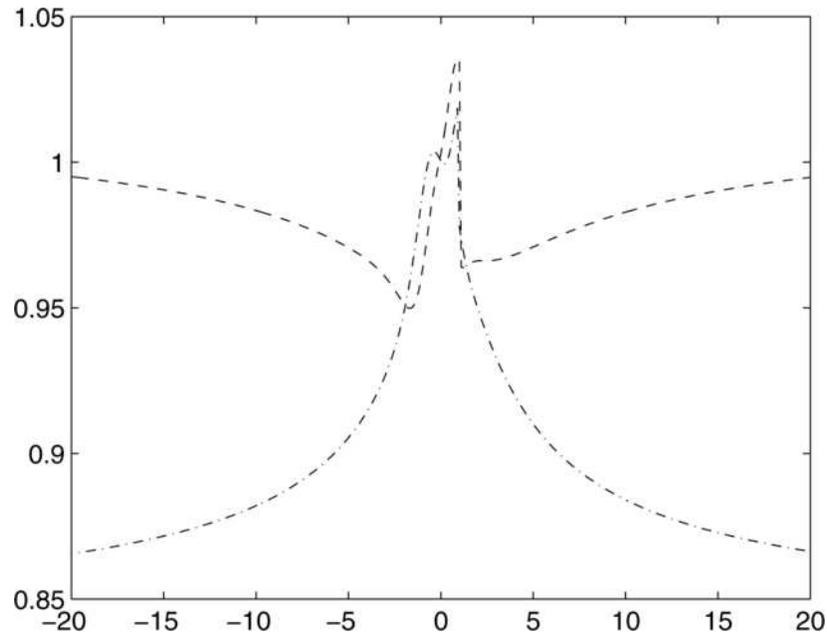

**Figure 3.** The tail error ratios described in (32) for $\hat{F}_R = \widehat{\Pr}_2$ (dashed) and $\hat{F}_R = \widehat{\Pr}_1$ (dash-dot).